\documentclass{amsart}
\usepackage{graphicx}
\usepackage{amsmath}
\usepackage{amsfonts}
\usepackage{amssymb}
\usepackage{latexsym}
\usepackage{array}
\usepackage{enumerate}
\usepackage{tikz}
\usepackage{floatflt}
\usepackage[T1]{fontenc}
\usepackage[default]{frcursive}
\newtheorem{theorem}{Theorem}

\theoremstyle{definition}
\newtheorem*{definition*}{Definition}
\theoremstyle{remark}
\newtheorem{remark}{\black{\bf Remark}}

\newenvironment{thm}%
{\color{black}\begingroup%
\begin{theorem}\vspace{-2pt}}
{\end{theorem}\endgroup}

\newenvironment{rem}%
{\color{black}\begingroup%
\begin{remark}}
{\end{remark}\endgroup}

\newenvironment{myproof}%
{\noindent\smallskip{\color{black}\proofname.}\begingroup%
}
{\endgroup}

\font\Times=cmr12 at 10truept

\newcommand{\Signe}{{\LARGE$\S$}}

\newcommand{\RL}{\mathcal{T}}

\newcommand{\Z}{\mathbb{Z}}
\newcommand{\mathL}{\mathcal L}

\newcommand{\N}{\mathbb{N}}

\renewcommand{\bar}{\overline}

\newcommand{\Id}{\textrm{Id}}

\newcommand{\QED}{\hfill\makebox[9pt]{\includegraphics[width=8pt]{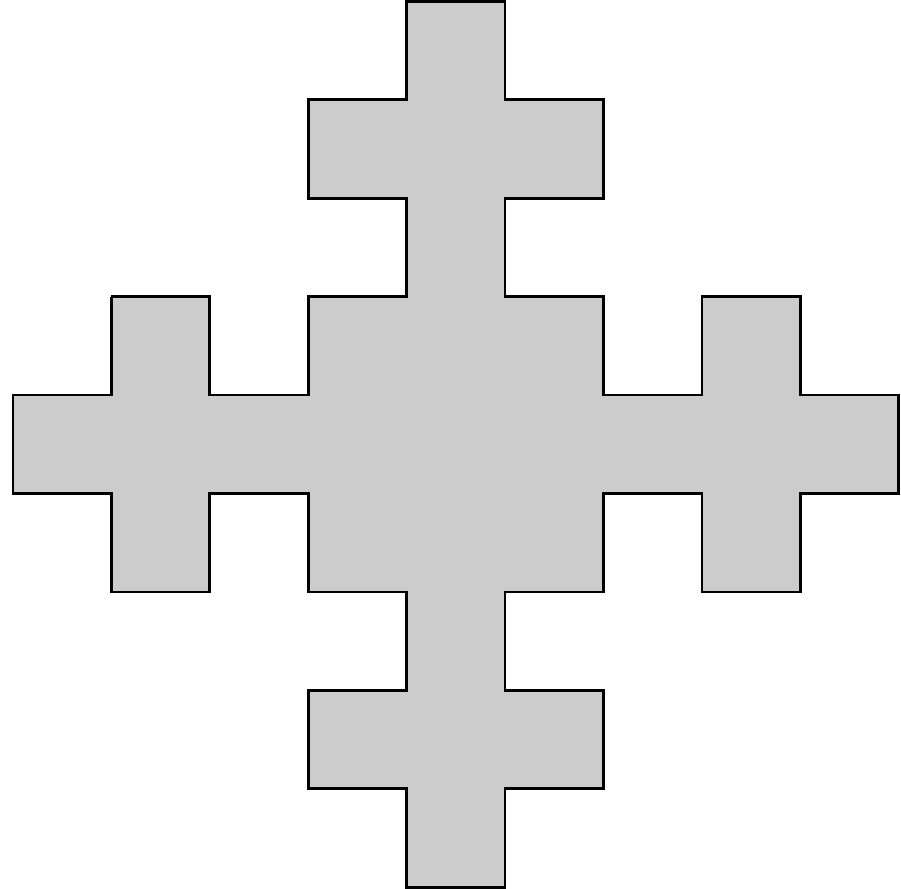}}}

\definecolor{NavyBlue}{rgb}{0.0, 0.0, .7}
\definecolor{citecolor}{cmyk}{0.2 0.2 0 0.7} 
\definecolor{itemcolor}{rgb}{0.1,0.4,0}
\newcommand{\red}[1]{{\color{red}#1}}

\newcommand{\black}[1]{{\color{black}#1}}

\newcommand{\mycite}[1]{{\color{citecolor}\cite{#1}}}
\newcommand{\myemail}[1]{{\color{black}\email{#1}}}
\newcommand{\myaddress}[1]{\noindent{\color{black}\address{#1}}}
\newcommand{\myparagraph}[1]{\stepcounter{section}\paragraph{\black{\Signe\bf\thesection.~ #1}}} 

\newcommand{\MMF}{M. Mend\`es France}
\newcommand{\SB}{S. Brlek}
\newcommand{\ABB}{A. Blondin Mass\'e}
\newcommand{\SL}{S. Labb\'e}

\newcounter{parcounter}

\begin{document}
\setcounter{parcounter}{0}

\renewcommand{\proofname}{\black{\bf Proof}}
\renewcommand{\abstractname}{\black{\scriptsize\bf Abstract}}

\title{ $\text{Complexity of the Fibonacci snowflake}$
}

  \author[\tt \tiny \ABB]{\tt \ABB}
\myaddress{\ABB, \SB, \SL, {\Times LaCIM}, Universit\'e du Qu\'ebec \`a Montr\'eal, C.P. 8888, Succ. Centre-ville,
   Montr{\'e}al (QC) Canada H3C 3P8}
 \myemail{\{blondin\_masse.alexandre, labbe.sebastien\}@courrier.uqam.ca,  \makebox[1cm]{      }\linebreak[4]   \makebox[2.8cm]{      }brlek.srecko@uqam.ca}
 \author[S.Brlek]{S. Brlek$^\dag$} 
\thanks{$^\dag$ 
With the support of {\Times NSERC} (Canada)}
  \author[S. Labb\'e]{S.Labb\'e}
 \author[M. Mend\`es France]{M. Mend\`es France}
 \myaddress{\MMF, D\'epartement de 
   math\'ematiques, UMR 5251, Universit\'e Bordeaux 1, 351 cours de la Lib\'eration, F-33405 Talence  cedex, France }
  \myemail{michel.mendes-france@math.u-bordeaux1.fr}

\dedicatory{A la m\'emoire de Benoit Mandelbrot}

\begin{abstract} The object under study is a particular closed curve on the square lattice $\Z^2$ related with the Fibonacci sequence $F_n$. It belongs to a class of curves whose length is $4F_{3n+1}$, and whose interiors by translation tile the plane. The limit object, when conveniently normalized, is a fractal line for which we compute first the fractal dimension, and then give a complexity measure.
\end{abstract}
\maketitle

{\color{NavyBlue}
\myparagraph{Introduction} 
 In a recent article \mycite{bblm} we defined and studied the properties of what we named a \emph{\color{red}Fibonacci snowflake}. Let us recall the main facts.
Consider the infinite  sequence $(q_n)_{n\in\N}$  defined recursively on the alphabet  $\RL=\{L=\text{left}, R=\text{right}\}$ by 
\begin{eqnarray*}\label{fibdef}
q_0&=&\varepsilon \quad\text{(the empty word}),\; q_1=R , \\
q_n &=& \begin{cases}
q_{n-1}q_{n-2}\;\; & \text{if $n\equiv2 \mod 3$,} \\
q_{n-1}\bar{q_{n-2}} & \text{if $n\equiv0,1 \mod 3$.} 
\end{cases}
\end{eqnarray*}
\noindent Here, $\bar{q_{n-2}}$ represents the word $q_{n-2}$ in which the letters $R$ and $L$ are permuted. The length $|q_n|$ of the word $q_n$ satisfies the Fibonacci relation
\[|q_n|=|q_{n-1}|+|q_{n-2}|\]
and therefore
\[|q_n|= \frac{1}{\sqrt{5}}\left(\frac{1+\sqrt{5}}{2}\right)^{n} -  \frac{1}{\sqrt{5}}\left(\frac{1-\sqrt{5}}{2}\right)^{n}.\]

To each word $w\in \RL^* \cup \RL^\N$ corresponds a polygonal line $\Pi$ on the lattice $\Z^2$. At the $n$-th vertex of $\Pi$, the $n$-th letter $w_n$ of $w$  indicates the direction of the next side. It may well happen that some sides are visited several times as for the case $w=L^4$ for example. The length of 
the polygonal line $\Pi$ corresponding to $w=L^4$ is : $|\Pi| =|L^4|+1.$ Actually, the following general property holds 
\begin{equation}\label{longueur}
|\Pi|=|w|+1, \; \forall w \in \RL^*.
\end{equation}

A finite polygonal line on $\Z^2$ is \emph{\color{red}closed} if both extremities coincide. The corresponding word is then said to be closed.  \emph{\color{red}Open} signifies non-closed; for example, $L^4$ is open while $L^3$ is closed. An open polygonal line is  \emph{\color{red} non-intersecting} if each of its vertices is attained once only. A \emph{\color{red} closed non-intersecting} polygonal line is one for which the only vertices visited twice are the extremities. Here are a few examples:
\begin{enumerate}[\color{itemcolor}(i)]
\item $L^4$ is neither closed nor non-intersecting;
\item $L^3$ is closed and non-intersecting;
\item $(LR)^n$ is open and non-intersecting;
\item $L^3RL^3$ is closed but not non-intersecting;
\item $(LRL)^3LR$ is closed and non-intersecting.
\end{enumerate}

\noindent
Given a word $w\not =\varepsilon$, we denote $w^{-\_}$ the word where the last letter is suppressed. For example, ${(LRL)^4}^{-\_}$ is closed and non-intersecting. 
In \mycite{bblm} we showed that the polygonal line $\Pi_n$ corresponding to the word ${(q_{3n+1})^4}^{-\_}$ is closed and non-intersecting (Theorem $10$) (see Figure \ref{F:basic fibo}). Therefore, the property \eqref{longueur} implies that 
\begin{equation}\label{lengthFS} |\Pi_n| = 4|q_{3n+1}|,  
\end{equation}
that is, the polygon $\Pi_n$ is composed of $4|q_{3n+1}|$ unit segments in $\Z^2$ . 
\begin{figure}[h!]
	\centering	
	\includegraphics[width=0.25\linewidth]{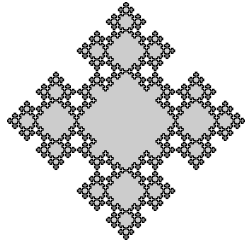}
	\caption{Fibonacci snowflake of order $n = 5$.}
	\label{F:basic fibo}
\end{figure}

\myparagraph{Fractal dimension}
The rather ``complicated'' associated polygon $\Pi_n$, when conveniently normalized, tends to a fractal line as $n\rightarrow\infty$ which we call \emph{\color{red}Fibonacci snowflake}. In our previous article we were aware of its complexity, yet we overlooked its fractal dimension even though we had all the information needed. We now prove:
\begin{thm} \label{fractaldim}
The fractal dimension of the Fibonacci snowflake is 
\begin{equation}
d= \frac{\log(2+\sqrt{5})}{\log(1+\sqrt{2})} = 1.637938210 \cdots
\end{equation}
\end{thm}
\begin{myproof}
First of all, we must see how to normalize the sequence of polygons $\Pi_n$ so that they stay in a bounded region  of the plane as $n\rightarrow\infty$. In our previous article we observed that the smallest square with sides parallel to the two axes $Ox,Oy$ and containing $\Pi_n$, has sides of size $2P(n+1)-1$ where the recurrence
\[P(0) = 0,P(1) = 1; P(n) = 2P(n-1) + P(n-2), \quad  \textrm{for $n > 1$},\]
defines the so-called  Pell numbers. Then they satisfy the relations
\begin{eqnarray*}
P(n+1)&=&\frac{2+\sqrt{2}}{4}  \left(1+\sqrt{2}\right)^n + \frac{2-\sqrt{2}}{4} \left(1-\sqrt{2}\right)^n \\
2P(n+1)-1&=&\frac{2+\sqrt{2}}{2}  \left(1+\sqrt{2}\right)^n + \frac{2-\sqrt{2}}{2} \left(1-\sqrt{2}\right)^n -1\\
&\sim &\frac{2+\sqrt{2}}{2}  \left(1+\sqrt{2}\right)^n .
\end{eqnarray*}

The polygon $\Pi_n$ is composed of $4|q_{3n+1}|$ unit segments and blows up as 
$n\rightarrow\infty$. But the \emph{\red{normalized}} polygon  $\frac{1}{2P(n+1)-1} \Pi_n$ stays bounded. It has $4|q_{3n+1}|$ sides each of length $(2P(n+1)-1)^{-1}$. The total $d$-dimensional normalized polygon has length
\[  \frac{4|q_{3n+1}|}{ (2P(n+1)-1)^{d}}\]
and therefore the fractal dimension of the Fibonacci snowflake is 
\[ d = \lim_{n\rightarrow\infty} \frac { \log(|q_{3n+1}|)}{\log(P(n+1))} =  \lim_{n\rightarrow\infty}   \frac { \log\left(\frac{1+\sqrt{5}}{2}\right)^{3n+1}}{\log(1+\sqrt{2})^n} =\frac{\log(2+\sqrt{5})}{\log(1+\sqrt{2})}. \quad\QED\] 
\end{myproof}

\myparagraph{A measure of complexity} There are obviously many ways to measure the complexity of a line. We mention below a measure related to the number of points of intersection of the figure with a \emph{random} straight line.

Recall a classical result due to Cauchy \mycite{cauchy}, Crofton \mycite{crofton}, Steinhaus \mycite{steinhaus}. Consider a plane rectifiable curve $\Gamma$ of length $|\Gamma|$ and whose convex hull $K$ has a frontier of length $|\partial K|$. A straight line $D$ is defined by 
\[x\cos\theta+y\sin\theta-\rho=0.\]  The probability measure considered is the uniform Lebesgue measure $d\rho.d\theta$ conditioned by the fact that $D$ intersects $\Gamma$. The result we alluded to is that the average number of intersecting points of $\Gamma$ with a random $D$ is 
\[N= \frac{2|\Gamma|}{|\partial K|}.\]

Applying this to $\Pi_n$, we see that the average number $N_n$ of intersection points of $\Pi_n$ with $D$ is 
\[N_n = 2\frac{2|q_{3n+1}|}{|\partial K_n|},\]
where $K_n$ is the convex hull of $\Pi_n$. Easily seen, $|\partial K_n|$ is of the order of $P(n+1)$. Therefore,

\begin{thm}\label{deux} We have
\[N_n \sim a \frac{(2+\sqrt{5})^n}{(1+\sqrt{2})^n} = a\left(\frac{2+\sqrt{5}}{1+\sqrt{2}}\right)^n .\]
\end{thm}

There is no difficulty to compute the constant $a=(1+\sqrt{5})$, but the important point is that the average number $N_n$ of intersection points increases exponentially, showing once more the high complexity of the Fibonacci snowflake. 
Theorem \ref{deux} could  obviously be considered as a corollary of the proof of Theorem \ref{fractaldim}. Both results are strongly related.\medskip

\begin{rem}\label{complexity} Let $\delta_n$ be the diameter of $\Pi_n$ (which is of the order of $P(n+1)$). Obviously, $2 \delta_n \leq |\partial K_n|\leq \pi\delta_n$. Theorem \ref{deux} states that the ratio $\mathL_n/\delta_n$ is of the order of $N_n$ where $\mathL_n$ is the length of $\Pi_n$. The difference 
\[\mathL_n-\delta_n=\delta_n\left(\frac{\mathL_n}{\delta_n}-1\right)\] measures the distance of $\Pi_n$ from being a straight line. The larger the ratio, more the curve meanders.
\end{rem}

\begin{rem} Let $F$ be a positive strictly increasing function. The ratio $F(\mathL_n)/F(\delta_n)$ is a measure of the complexity of $\Pi_n$. When $F=\log$ we obtain Theorem \ref{fractaldim} and when $F=\Id$ we obtain Theorem \ref{deux}. The choice $F=\exp$ leads to $\exp(\mathL_n-\delta_n)$ i.e. Remark \ref{complexity}.
\end{rem}
\bigskip
\myparagraph{Entropy}
Let $p_j^{(n)}$ be the probability that a random straight line intersects  $\Pi_n$ in exactly $j$ points, given that the line meets $\Pi_n$. The associated entropy    $h_n$ is by  definition   
\[h_n =-\sum p_j{(n)}  \log (p_j{(n)}).\]
In \mycite{Mmf} it is shown that
\[h_n\leq \log \left(\frac{2\mathL_n}{|\partial K_n|} \right)+\left(1- \frac{|\partial K_n|}{2\mathL_n}\right) \log\left(\frac{2\mathL_n}{2\mathL_n-|\partial K_n|} \right). \]
The second term in the righthand side is positive, less than $1$ and tends to $0$ with $\frac{|\partial K_n|}{\mathL_n}.$

Since $\mathL_n$ is of the order of $(2+\sqrt{5})^n$ and $|\partial K_n|$ is of the order of $(1+\sqrt{2})^n$, we have 
\[ h_n\leq n \log \left( \frac{2+\sqrt{5}}{1+\sqrt{2}}\right) + {\mathcal O} (1).\]
We have thus established an upper bound for the complexity of $\Pi_n$: 
\begin{thm}  
\[ \limsup_{n \to \infty} \frac{1} {n}h_n  \leq \log \left( \frac{2+\sqrt{5}}{1+\sqrt{2}}\right) .\]
\end{thm}


\renewcommand{\bibname}{\black{\bf References}}

}
\end{document}